\documentclass[twocolumn]{article}

\usepackage{color}

\usepackage{amsmath}

\begin{document}

\title{How Could Polyhedral Theory Harness Deep Learning?\footnote{Presented at the Scientific Machine Learning Workshop, U.S. Department of Energy
Office of Advanced Scientific Computing Research, January 30 -- February 1, 2018}}
\author{Thiago Serra ~ Christian Tjandraatmadja \\  {\small Carnegie Mellon University} \and Srikumar Ramalingam \\ {\small The University of Utah}}
\date{}

\maketitle

\subsection*{Abstract}
The holy grail of deep learning is to come up with an automatic method to design optimal architectures for different applications. In other words, how can we effectively dimension and organize neurons along the network layers based on the computational resources, input size, and amount of training data? We outline promising research directions based on polyhedral theory and mixed-integer representability that may offer an analytical approach to this question, in contrast to the empirical techniques often employed.

\subsection*{Introduction}
Deep neural networks can be seen as function approximators, and a network is considered to be better at modeling certain problems if it has better expressiveness. One hypothesis for the expressiveness of these networks is that more layers can help breaking down the input space into exponentially more regions that are evaluated differently~\cite{Bengio2009}. Those regions are often linear 
in practice because some of the most useful and popular activation functions such as rectifier linear units (ReLU) and maxouts are piecewise-linear functions~\cite{maxoutnetworks13,Nair2010}. 
This led to a novel stream of work 
fueled by integer programming
and its related polyhedral studies, a contrasting area where some theory waited decades for computers that could exploit it~\cite{GMI}.

Those linear regions can be represented  
as the union of polyhedra consisting of all combinations of active hyperplanes for each neuron, 
from which we discount infeasible terms and proper subsets. 
Counting the remaining terms and understanding their impact on effectiveness is a flourishing area of research. 

\subsection*{Research Directions}

\subsubsection*{1. Bounds on the linear regions}

Current bounds on linear regions are based on the theory of hyperplane arrangements~\cite{Zaslavsky1975} and on constructions increasing the realizable number of regions. Early results have validated that deeper networks with uniform width across layers can be more expressive in terms of the number of linear regions for small input dimensions~\cite{Montufar2014,Pascanu2013}. Subsequent work has improved the asymptotic bounds for uniform width~\cite{Raghu2017} and total number of neurons~\cite{Arora2016}, and later refined these figures~\cite{Montufar2017,BoundingCounting}. A closer look at the polyhedra defined by the input space and transformed along the layers has been central in recent work. For example, by using zonotopes to construct networks with more regions~\cite{Arora2016} and exploiting the dimension of the space effectively discriminated by the network to bound the number of hyperplane arrangements across all regions~\cite{Montufar2017} and on each one of them~\cite{BoundingCounting}. 
We believe that further work could validate the insights from those bounds to design more expressive networks.  
For example, the current bounds suggest that there are particular depths that maximize the number of linear regions for a given input size and number of neurons.

\subsubsection*{2. Enumeration of linear regions}

The collection of linear regions can be modeled with mixed-integer linear formulations, where each linear region defines a different assignment of integer variables. These formulations have been used to generate adversarial examples~\cite{Cheng2017,Fischetti2017} as well as to count linear regions for bounded inputs~\cite{BoundingCounting}. 
While the exact enumeration of linear regions is currently impractical for the large deep networks that are often employed in applications, 
the connection with discrete optimization methods may lead to better techniques in the future. 
Another potential direction would be to develop upper and lower approximations that scale better in practice. 
The progress in counting can provide novel insights in improving the bounds and vice versa. 

\subsubsection*{3. What does the number of regions tells us?}

So far, there is only anecdotal evidence connecting potential expressiveness to how effective these networks are in practice. 
Extensive experiments could aim for statistical significance between the number of linear regions and the accuracy of the  network.

This line of research can be particularly challenging because expressiveness may also relate to the potential for over-fitting. We have previously hypothesized that the network training is not likely to generalize well if there are so many regions that each point can be singled out in a different region, in particular if regions with similar labels are unlikely to be compositionally related. If successful, this line of work could help defining network configurations that would be better suited for the available amount of data.

\subsubsection*{4. The shape of linear regions}

A different take on expressiveness would be to study the shape of linear regions. This topic has been only lightly approached in the literature, but with interesting insights concerning the compositional benefits of having more than one layer~\cite{Raghu2017,BoundingCounting}. By understanding if and to what extent the input space can be transformed along additional layers, it might be possible to determine the necessary depth to outline the input space in certain ways, thereby assessing the computational resources needed for particular tasks. 

Furthermore, the existence of mixed-integer formulations imply that ReLU and maxout activations are mixed-integer representable~\cite{MIR}. Therefore, although neural networks are universal approximators~\cite{Cybenko1989}, some negative results can be derived when the total number of neurons is finite.

\subsection*{Conclusion}

We believe that further work on those topics will result in actionable insights that could orient future practice in a more systematic way. 
This could help practitioners to understand the impact of their design decisions on different applications.

{\footnotesize
\bibliography{DeepLearningExp}
\bibliographystyle{plain}}

\end{document}